# NUMERICAL ANALYSIS OF MINIMAL BETA-SEQUENCES ASSOCIATED WITH A FAMILY OF ENTIRE FUNCTIONS

ALLAN M. DIN [*)] and LORENZO MONETA [**)]

ABSTRACT. The Riemann Xi-function $\Xi(t) = x(1/2 + it)$ is a particularly interesting member of a broad family of entire functions which can be expanded in terms of symmetrized Pochhammer polynomials depending on a certain scaling parameter $b$. An entire function in this family can be expressed as a specific integral transform of a function A(x) to which can be associated a unique minimal beta-sequence $b_{\min,n} \to \infty$ as $n \to \infty$, having the property that the Pochhammer polynomial approximants $\Xi_n(t, b_n)$ of order n to the function $\Xi(t)$ have real roots only in t for all n and for all $b_n \geq b_{\min,n}$.

The importance of the minimal beta-sequence is related to the fact that its asymptotic properties may, by virtue of the Hurwitz theorem of complex analysis, allow for making inferences about the zeros of the limit function $\Xi(t)$ in case the approximants $\Xi_n(t, b_n)$ converge. The objective of the paper is to investigate numerically the properties, in particular the very large n properties, of the minimal beta-sequences for different choices of the function A(x) of compact support and of exponential decrease, e.g. like the Riemann case.

___________

# 1. Introduction

Entire functions are complex-valued functions which are holomorphic in the whole complex plane and they can be expanded as a power series which converges uniformly on compact sets [see e.g. Ref 9]. Such functions may be represented by a product involving its zeroes and there exists a vast literature about the detailed characteristics of these zeros [see e.g. Ref 8]. However, one of the most well-known entire functions, the Riemann Xi-function $\Xi(t) = x(1/2+it)$ [see Ref. 13], has for more than 150 years resisted all attempts towards proving the Riemann Hypothesis stating that all non-trivial zeros of $\Xi(t)$ are located on the critical line $s = 1/2+it$ with t real.

It happens that the $\Xi$-function is a member of a rather broad family of entire functions, as defined by a certain integral representation, and some of these functions have quite a number of properties in common with the Riemann $\Xi$-function. This is particularly the case for a certain Bessel function of the third kind, and since this function only has zeros located on the critical line, one might easily be led to speculate why this would not be the case also for the Riemann $\Xi$-function. Unfortunately there are also examples of family members which have other similarities with the Riemann $\Xi$-function, but which nevertheless have complex zeros off the critical line.

There exists a great number of different attempts to prove the Riemann Hypothesis and very often they result in the statement of another equivalent hypothesis, frequently concerning a field of mathematics more or less unrelated to complex analysis [for a review see Ref. 3]. Concerning approaches using complex analysis, it is somewhat surprising that historically there has been no serious attempt to investigate the validity of the Riemann Hypothesis in terms of the well-known Hurtwitz theorem of complex analysis [see Ref.9], which in a quite manifestly concerns the issue at hand about zeros of analytic functions. This theorem states that if an analytic function is a limit of a sequence of analytic functions, uniformly convergent on a compact subset of the complex plane, then any zero of the function in the subset must be a limit of the zeros of the sequence functions. In particular, if the zeros of the approximants were all real, then the same would be true for the limit function.

One explanation for the apparent lack of interest in trying to apply the Hurtwitz theorem may be related to the fact that, apart from knowing that a certain entire function may theoretically be represented by a product involving its zeroes, then it is in general not obvious how to construct explicitly an appropriate sequence of analytic functions converging to the given function. For the problem of the Riemann Hypothesis, of course "appropriate" here means that the approximating functions should have only real zeros so as to be able to infer that the same would be the case for the limit function.

In some simple cases, like for the function sin(t), the situation is pretty clear since we know explicitly the infinite product representation for the entire function

$$\sin(t) = t \prod_{k=1}^{\infty} (1 - \frac{t^2}{k^2 p^2}) \qquad (1.1)$$

Here we can easily define a sequence of finite products, i.e. polynomials in t, with real zeros which converges to sin(t), and this function can therefore only have real zeros, a fact which is of course already known.

In general, it appears to be more feasible to represent an entire function directly as some infinite sum of simple functions like powers or polynomials. For the above example we have of course the standard Taylor expansion

$$\sin(t) = \sum_{k=0}^{\infty} (-1)^k \frac{t^{2k+1}}{(2k+1)!} \qquad (1.2)$$

but here the application of the Hurwitz theorem to the sequence of polynomial approximants does not provide any useful information about real zeros. In fact the partial sums of the above series have mostly complex zeros and there seems to be no obvious generic way of correcting this situation.

The Hurwitz theorem might nevertheless be deemed to be relevant for the stated purpose if it were somehow possible to define a sequence of approximants depending on a certain parameter which could be suitably adjusted to the make the approximants have only real zeros. Surprisingly, this is actually feasible for the family of entire functions mentioned above, including the Riemann $\Xi$-function, because it happens that these functions can be expanded in terms of a uniformly convergent series of symmetrized Pochhammer polynomials depending on a certain scaling parameter $b$ [see Ref 4]. Moreover, for large values of $b$ it appears that these polynomial approximants have real roots only.

The objective of the present paper is precisely to undertake a numerical investigation of increasing beta-sequences $b_n$ which are such that a bona fide polynomial approximant $\Xi_n(t, b_n)$ to a given function $\Xi(t)$ has only real roots. Specifically it is of interest to investigate beta-sequences which grow as slowly as possible, so-called minimal beta-sequences. We will show how appropriate root finding algorithms can be used to calculate numerically such minimal beta-sequences for large n, as well as to find their associated t-root sequences which appear to have some very peculiar properties. Moreover we will investigate numerically certain functional indicators which provide quite valuable insights into the underlying asymptotic mechanism for the observed behaviour of these sequences.

However, we will not here discuss the more formal asymptotic analysis of these sequences since this is done elsewhere [see Ref. 4]. Such a formal analysis is warranted for several reasons. Firstly, there is the issue of analyzing whether the numerical results for the large n behaviour of the sequences correspond to the true asymptotic behaviour, insofar as this can be disentangled analytically. Secondly, only formal analysis can disentangle under which conditions the bona fide approximants actually converge to the limit function. Finally, the open question to be settled analytically is whether all conditions are actually united so as to be able to apply the Hurwitz theorem to the case of the Riemann $\Xi$-function, and thus eventually confirm the validity of the Riemann Hypothesis based entirely on complex analysis.

**2. The family of entire function**

The Riemann xi-function $x(s)$ is related to the Riemann zeta-function $V(s)$ by [see e.g. Ref. 13]

$$x(s) = \Gamma(s/2+1)(s-1)p^{-s/2}V(s) \qquad (2.1)$$

It is an entire function of s which has the same (non-trivial) zeros as $V(s)$ in the critical strip $0<\text{Re}(s)<1$ and it fulfils the functional equation $x(s) = x(1-s)$. One has the explicit representation

$$x(s) = \int_1^\infty dx A(x)(x^{-s/2} + x^{(s-1)/2}) \tag{2.2}$$

where A(x) is given in terms of the elliptic theta function

$$y(x) = (J_3(0, \exp(-px)) - 1)/2 = \sum_{n=1}^\infty \exp(-n^2 px) \tag{2.3}$$

by

$$A(x) = 2\frac{d}{dx}(x^{3/2} y'(x)) = \sum_{n=1}^\infty (2n^4 p^2 x - 3n^2 p) x^{1/2} \exp(-n^2 px) \tag{2.4}$$

One notes for this A(x), defined a priori on the interval $[1, \infty]$, that we have A(x)>0 and that it is bounded by a power of x times $\exp(-px)$. We will here be interested in studying the broad family of entire functions defined by the above integral representation in terms of a certain function A(x) which is just required to be non-negative function for all $x \geq 1$ with A(1)>0 and decreasing exponentially. i.e. like $\exp(-ax^b)$ with a>0 and b>0, or faster for large x. Actually there is no loss of generality by only studying A(x) of exponential order one, i.e. b=1 [see Ref. 4]. For simplicity, we will continue using the notation $x(s)$ and refer to it as the Riemann $x(s)$ when the specific Riemann representation of A(x) is used.

We recall that the Riemann A(x) fulfills the inversion transformation relation

$$A(x) = x^{-3/2} A(1/x) \tag{2.5}$$

which means that the combination $A_I(x)$ defined by

$$A_I(x) = x^{3/4} A(x) \tag{2.6}$$

is invariant under the operation $x \to 1/x$. This property is known to be instrumental for the fact that the Riemann $x(s)$ has infinitely many zeros on the critical line, but we will also investigate some A(x) cases for which this inversion invariance is not fulfilled.

To make the underlying symmetry in the complex plane manifest, we will make the usual change of parameters $s = 1/2 + it$ so that the critical line $\text{Re}(s) = 1/2$ is the line t = real and the critical strip $0<\text{Re}(s)<1$ is $|\text{Im}(t)| < 1/2$. We also introduce the Xi-function $\Xi(t) = x(1/2 + it)$ so that the integral representation becomes

$$\Xi(t) = \int_1^\infty dx A_I(x) x^{-1} (x^{it/2} + x^{-it/2}) \tag{2.7}$$

The basic character of the mentioned inversion symmetry become clearer when the integral representation of $\Xi(t)$ is recast in the standard form of a cosine transformation after the change of variable $x = e^{2y}$:

$$\Xi(t) = 4\int_0^\infty dy A_I(e^{2y})\cos(ty) \qquad (2.8)$$

The 1/x symmetry can now be understood as a reflection symmetry $y \to -y$ of $A_I(e^{2y})$. Generally if $A(x)$ is chosen so that $A_I(x)$ is invariant under $x \to 1/x$, then the resulting $\Xi(t)$ is likely to have an infinite number of real zeros. This is so, for example, if $A_I(x)$ depends simply on the invariant combination $x + 1/x$, say like $\exp[-a(x+1/x)]$ which leads to a $\Xi(t)$ involving Bessel functions of the third kind, as considered historically when examining certain bona fide approximations to the Riemann $\Xi(t)$ function [e.g. Ref. 13].

If we choose a=1 then we get the following simple and representative member of the A(x) family:

$$\Xi(t) = 2K_{it/2}(2) \text{ for } A_I(x) = Exp\left[-(x+\frac{1}{x})\right] \qquad (2.9)$$

Different Bessel function examples of this type, e.g. for a<>1 and for power factors multiplying the exponential term in $A_I(x)$, appear to have not only an infinite number of real only zeros but they also have an asymptotic density distribution similar to the one of the Riemann $\Xi(t)$. If we break inversion symmetry by choosing simply $A_I(x) = \exp(-ax)$, with a=1 then we recover another interesting family member:

$$\Xi(t) = \Gamma(it/2,1) + \Gamma(-it/2,1) \text{ for } A_I(x) = Exp[-x] \qquad (2.10)$$

which is a symmetrized incomplete gamma function which has a real zero at infinity only. It is also worth noticing that if one only retains the first term in the theta function expansion of the Riemann A(x), then the resulting $\Xi(t)$ has only one real zero. As more terms are retained, more and more real zeros appear. The corresponding $\Xi(t)$ function simply corresponds to a sum of symmetrized incomplete gamma functions.

Another instructive example appears if $A_I(x)$ is chosen simply to be 1 with compact support in the x-interval $[1, \exp(2w)]$ for any $w>0$, with $w=1$ as a good representative example. Then we find

$$\Xi(t) = 4\sin(wt)/t \text{ for } A_I(x) = 1, x \in \left[1, e^{2w}\right] \qquad (2.11)$$

If we choose $A_I(x)$ to be one at x=1 and decreasing linearly to zero on the same interval as above, then for $w=1$ we find $\Xi(t) = 8(\sin(t/2)/t)^2$ and we have a limiting case when the real zeros (non-simple) are about to disappear completely. Another interesting intermediate case arises when $A_I(x)$ is chosen to be

cos(log(x)$p$/4) on the same interval as above for $w=1$, which leads to a quadratically decreasing $\Xi(t) = 2p \cos(t)/((p/2)^2-t^2)$ having an infinite number of simple zeros.

We will also investigate some more "exotic" members of the A(x) function family which are cleverly engineered examples of Dirichlet series, known in the Riemann zeta-function literature for quite some time, which have properties comparable to the Riemann series but which nevertheless explicitly exhibit complex zeros of their corresponding $\Xi(t)$ function. The examples to be studied below are two different one-parameter families of $A_I(x)$ with this property. The first one, denoted $A^{(1)}$, is related to the Ramanujan tau function [Ref. 2] and is explicitly given by

$$A^{(1)}(x,k) = \left[ x^{1/8} \exp(-p\sqrt{x}/12) \prod_{n=1}^{\infty} (1-\exp(-2pn\sqrt{x})) \right]^k \tag{2.12}$$

where k is a positive integer different from 1, 2, 3, 4, 6, 8, 12 and 24. Such an A function is of exponential order 1/2, decreases monotonously like the Riemann A(x), and corresponds to a $\Xi(t)$ function which has infinitely many real zeros, but which also has complex zeros outside the critical strip.

A second one-parameter A(x) family of exponential order 2, denoted $A^{(2)}$, is given by [Ref. 1]:

$$A^{(2)}(x,k) = \sum_{n=1}^{\infty} \left[ kA^{(2,1)}(x,n) - c(n)A^{(2,2)}(x,n) \right]$$
$$A^{(2,1)}(x,n) = (4p^2n^4x^4 - 6pn^2x^2)\exp(-pn^2x^2) \tag{2.13}$$
$$A^{(2,2)}(x,n) = (4p^2n^4x^4/25 - 6pn^2x^2/5)\exp(-pn^2x^2/5)$$

where $c(n)$ is the quadratic Dirichlet character modulo 5 ($c(1) = c(4) = 1, c(2) = c(3) = -1, c(5) = 0$) and $k \geq 4$. For k>4 the corresponding $\Xi(t)$ function has no real zeros and presumably an infinite number of complex zeros inside and outside the critical strip.

## 3. The Pochhammer polynomial expansion

The Pochhammer polynomials $P_k(s)$ of degree k are defined by

$$P_k(s) = \prod_{j=1}^{k} (1 - \frac{s}{j}) = \frac{\Gamma(k+1-s)}{\Gamma(k+1)\Gamma(1-s)} \tag{3.1}$$

with $P_0(s) = 1$, $P_1(s) = 1 - s$, and $P_2(s) = 1 - 3s/2 + s^2/2$, etc. The Pochhammer polynomials have a simple generating function

$$(1-e)^s = \sum_{k=0}^{\infty} P_k(s+1)e^k \tag{3.2}$$

with the series being absolutely convergent for $|e|<1$.

In general, the $\Xi(t)$ function has the convergent expansion in Pochhammer polynomials depending on an arbitrary positive parameter $b$ [see Ref. 4]

$$\Xi(t) = \sum_{k=0}^{\infty} b_k(b) P_k^+(t/b) \quad (3.3)$$

where we have explicitly introduced the symmetrized Pochhammer polynomial

$$P_k^+(t) = (P_k(it) + P_k(-it))/2 \quad (3.4)$$

which is just the even part of $P_k(it)$, and the coefficient $b_k(b)$ is given by the integral

$$b_k(b) = 2\int_1^{\infty} dx A_I(x) x^{-1} x^{-b/2} (1-x^{-b/2})^k = \frac{4}{b} \int_0^1 dy A_I(y^{-2/b})(1-y)^k \quad (3.5)$$

Let us note that we have the general decomposition

$$P_k(it) = P_k^+(t) + i P_k^-(t) \quad (3.6)$$

where $P_k^+(t)$ is the even part of $P_k(it)$, as defined above, and $P_k^-(t)$ is the odd part of $P_k(it)$ divided by i. For reference let us write out explicitly the expressions for the lowest k-values:

$$P_0^+(t) = 1, P_1^+(t) = 1, P_2^+(t) = 1 - t^2/2, P_3^+(t) = 1 - t^2 \quad (3.7)$$

$$P_0^-(t) = 0, P_1^-(t) = -t, P_2^-(t) = -3t/2, P_3^-(t) = -t(11-t^2)/6 \quad (3.8)$$

It is quite easy to convince oneself that $P_k^+(t)$ and $P_k^-(t)$ have distinct real roots only.

Starting from the convergent expansion of $\Xi(t)$ in symmetrized Pochhammer polynomials, valid for any $b>0$, we can simply define the polynomial approximants by

$$\Xi_n(t,b) = \sum_{k=0}^{n} b_k(b) P_k^+(t/b) \quad (3.9)$$

which have the property that $\Xi_{2n}(t,b)$ as well as $\Xi_{2n+1}(t,b)$ are even alternating polynomials of degree n in $t^2$, as it can easily be seen from the general properties of the polynomials $P_k^+(t)$. For fixed $b$ these approximants converge to $\Xi(t)$, but even if the symmetrized Pochhammer polynomials have real roots only, we would not a priori expect the approximants to have real roots, keeping in mind the simple example of the introduction.

Surprisingly it appears that for low n, then the approximant $\Xi_n(t,b)$ has mostly real roots, in contrast to the standard situation of partial sums of Taylor expansions. But starting from a certain n, then complex roots of $\Xi_n(t,b)$ emerge and remain for higher n. However, it can be shown [see Ref. 5] that for any given n, then it is always possible to find a $b$ so that $\Xi_n(t,b)$ has only real roots from this $b$ and onwards. Because of this very special feature, one may infer that there exists increasing beta-sequences $b_n$ such that $\Xi_n(t,b_n)$ has real roots only for all n. Moreover, one can assert that for a given A(x), and thus $\Xi(t)$, then there exists a unique minimal beta-sequence $b_{\min,n}$ which is simply defined as an infimum over the above beta-sequences. This minimal beta-sequence fulfills a specific difference equation which will be described in the next section, and which will be studied numerically in the subsequent sections.

**4. Properties of minimal beta-sequences**

It is not too complicated to find the basic properties of minimal beta-sequences when keeping in mind how they emerge as a limiting situation of real roots only of a polynomial depending on the $b$ parameter. Simply, if a polynomial has distinct real roots for a certain (sufficiently large) $b$, then as the real parameter $b$ is decreased, the polynomial roots must follow continuous trajectories in the complex plane [see Ref. 5, 11 for a review of polynomial root dynamics]. At a certain point, one will generally see two real roots coalesce and then split into a pair of complex roots. The minimal beta-sequence case therefore corresponds precisely to a real double-root situation of the polynomial approximant.

Let us therefore suppose that for a given n, we have already found the values $b_n$ and $t_n$ corresponding to the situation where the polynomial $\Xi_n(t,b)$ has a real double root, such that for all $b > b_n$ the roots are all real and distinct while for $b < b_n$ complex roots start to appear. Consequently we have that both $\Xi_n(t,b)$ and its first derivative with respect to t, denoted by $\Xi_n'(t,b)$, are both zero at this point:

$$\Xi_n(t_n, b_n) = 0 \text{ and } \Xi_n'(t_n, b_n) = 0 \qquad (4.1)$$

If we now increment n by one, we have

$$\Xi_{n+1}(t,b) = \Xi_n(t,b) + b_{n+1}(b) P_{n+1}^+(t/b) \qquad (4.2)$$

and the values of $b_n$ and $t_n$ will also have to be incremented

$$b_{n+1} = b_n + \Delta b_n, \, t_{n+1} = t_n + \Delta t_n \qquad (4.3)$$

To first order in the increments one finds the following difference equations, or recursive equations, which we will refer to below as the beta and t-increment equations:

$$\Delta b_n = -\frac{\Xi_{n+1}(t_n, b_n)}{\frac{\partial}{\partial b}\Xi_{n+1}(t_n, b_n)} \quad (4.4)$$

$$\Delta t_n = -\frac{\Xi'_{n+1}(t_n, b_{n+1})}{\Xi''_{n+1}(t_n, b_n)} \quad (4.5)$$

These equations simply encompass the classical Newton root finding method, and if they are iterated in the indicated order, keeping the polynomial order fixed, they provide a straightforward numerical approach to finding the double root values at any desired precision. The first order approximation is typically correct to within 3% and the Newton type iteration therefore converges rather quickly. For low n in the range up to n=400 it is quite feasible to do computer calculations using standard PC software like Wolfram Mathematica, but for much higher n of the order of one million it is necessary to adapt the approach to a mainframe computer platform. The results of these two approaches to the numerical analysis are presented in the following sections.

For rather low n, say in the range 4-50, it is quite feasible to numerically calculate all roots of the polynomial $\Xi_n(t,b)$, of course after first calculating the coefficients $b_n(b)$ by numeric integration and the symmetrized Pochhammer polynomials by recursion. When doing this, it is straightforward to observe explicitly that $\Xi_n(t,b)$ for large $b$ has only distinct real roots and that, when $b$ is decreased, two real roots will be getting closer and closer, and will finally coalesce to become a real double root at $t_n$, before becoming a pair of complex roots. For low n, this fusion of two real roots to a double root is seen to happen in the lowest root range, but as n increases then the fusion point $t_n$ will be seen to move higher and higher.

For large n, say of the order of thousands or millions, we use a mainframe computing platform and apply a multi-dimensional numerical root finder algorithm to solve iteratively the two double root equations (4.1). The algorithms provided by the ROOT framework [see Ref.12] are used to perform the numerical calculations and the beta-coefficients required to evaluate the Xi-function are obtained by integrating numerical the A(x) function using an adaptive integration algorithm implemented in ROOT using the GNU Scientific Library [see Ref.7].

The solution of the two equations is found using a multi-dimensional root finding algorithm based on a hybrid Newton method [see Ref.10] in the implementation provided by the GNU Scientific Library. The algorithm requires an initial guess for the t and beta values, which for low n may be found by explicit calculation of all roots, and subsequently it uses the solution found at the previous iteration. However, when t-jumps are encountered (see Section 5), a special procedure must be applied whereby one shifts the initial parameter guess to higher t-values. In practice, a t-jump is signaled by observing a failure of the algorithm in finding a solution within the specified tolerance, and then the iterative procedure of shifting t to higher values continues until the algorithm converges. As an alternative to the iterative t-jump search, it is possible, for a given A(x) function, to anticipate quite precisely the location of the t-jump using the known expression of the corresponding Xi-function.

The application of the above implementation of the Newton method allows for pushing the computation up to values of n of the order of a million, but eventually the numerical evaluation of the b-coefficient integral becomes very CPU time consuming (of the order of days). In particular, the double numerical accuracy implies a heavy computing time penalty, so to push the calculations to even higher n-values, more efficient and precise methods for numerical integral evaluation would have to be used.

Before entering in the details of the numerical results for the beta and t-sequences, it is worth noticing that when n becomes very large then the polynomial increments are very small and therefore the first order approximation becomes better and better. Consequently, the formal asymptotic analysis of the beta and t-sequence can simply be applied using the above equations without iteration, in the following concise form explicitly implementing the double root conditions $\Xi_n(t_n, b_n) = 0$ and $\Xi'_n(t_n, b_n) = 0$:

$$\Delta b_n = -b_{n+1}(b_n) Q_{n+1}(t_n, b_n) \qquad (4.6)$$

$$\Delta t_n = -\frac{b_{n+1}(b_n) P^+_{n+1}(t_n/b_n)}{\Xi''_{n+1}(t_n, b_n)} \frac{\partial}{\partial t} \log(Q_{n+1}(t_n, b_n)) \qquad (4.7)$$

$$Q_n(t_n, b_n) \equiv \frac{P^+_{n+1}(t_n/b_n)}{\frac{\partial}{\partial b}\Xi_{n+1}(t_n, b_n)} \qquad (4.8)$$

Here we have introduced explicitly the ratio $Q_{n+1}(t_n, b_n)$ for which the interplay of the numerator and the denominator is instrumental for the large n asymptotic analysis. There are other simplifications of the asymptotic region which are important. Firstly, the symmetrized Pochhammer polynomial $P^+_{n+1}(u_n)$, where have put $u_n = t_n/b_n$, can be approximated by using the asymptotics of gamma function to write

$$P_n(it) \approx n^{-it}/\Gamma(1-it) \qquad (4.9)$$

Consequently

$$P^+_n(u) \approx \frac{\cos(u \log(n) + j(u))}{r(u)} \qquad (4.10)$$

where we have set

$$r(u) = Abs(\Gamma(1-iu)) = \left(\frac{pu}{\sinh(pu)}\right)^{1/2}, \quad j(u) = Arg(\Gamma(1-iu)) \approx Cu \qquad (4.11)$$

and the last relation is valid for small u with C being the Euler constant. Thus $P^+_{n+1}(u_n)$ behaves asymptotically in an oscillatory manner depending on $u_n$. Secondly, the beta-derivative term appearing in the asymptotic beta-increment equation can be seen [see Ref. 4] to be related asymptotically to the first derivative of the $\Xi(t)$ function in the following way:

$$b_n \frac{\partial}{\partial b} \Xi_{n+1}(t_n, b_n) \approx \Xi^*(t_n + t_n) \tag{4.12}$$

where $t_n$ is a t-jump quantity to be discussed below, and

$$\Xi^*(t) \equiv \Xi(t) + t\Xi'(t) \tag{4.13}$$

After calculating the $b_n$ and $t_n$ sequences numerically by the iterative procedure, it will be useful to exhibit explicitly the behaviour of the above terms so as to better understand analytically the underlying mechanism of the asymptotic increment equations.

## 5. Reference cases with compact support

For the case of $A_1(x) = 1$ on the interval $[1,\exp(2)]$, then we have $\Xi(t) = 4\sin(t)/t$. Moreover the coefficients $b_k(b)$ in the polynomial expansion of $\Xi_n(t, b)$ can be evaluated exactly as

$$b_k(b) = \frac{4}{b(k+1)}(1 - e^{-b})^{k+1} \tag{5.1}$$

so the numerical calculations are somewhat simpler and faster to implement in this case.

We may start the evaluation of the sequences with the lowest possible n=4 and find first $b_4 \sim 0.074$ and $t_4 \sim 3.462$, and then $b_5 \sim 0.268$ and $t_5 \sim 3.471$. This represents a small incremental change of t and therefore the t-increment equation is directly applicable without any further consideration. But for n=6, then one notes that t jumps to a significantly higher value $t_6 \sim 5.15$, while beta drops to $b_6 \sim 0.126$. This phenomenon of t-jump is of course also quite manifest if one would calculate all roots of the polynomial for n=6, and then check which roots coalesce as beta is decreased from a certain high value.

From the point of view of a computerized procedure, it is therefore clear that the local Newton like iteration method at n = 6 must be extended so as to be able to find a new polynomial extremum point, which may be situated at a certain distance from the preceding double root. However, it is quite well-known how to implement an appropriate non-local extension of the method. Subsequently one may proceed with the simplified iterative procedure and note that both β and t increase smoothly until n = 17 where there is a new t-jump from $t_{16} \sim 5.69$ to $t_{17} \sim 7.93$, while β drops very slightly from $b_{16} \sim 1.22$ to $b_{17} \sim 1.19$. This process continues with t-jumps at n = 40, 87, 185, 386 and further on. The graphs for $b_n$, $t_n$ and $u_n = t_n / b_n$ in the n-range up to 400 is presented in Fig. 5.1.

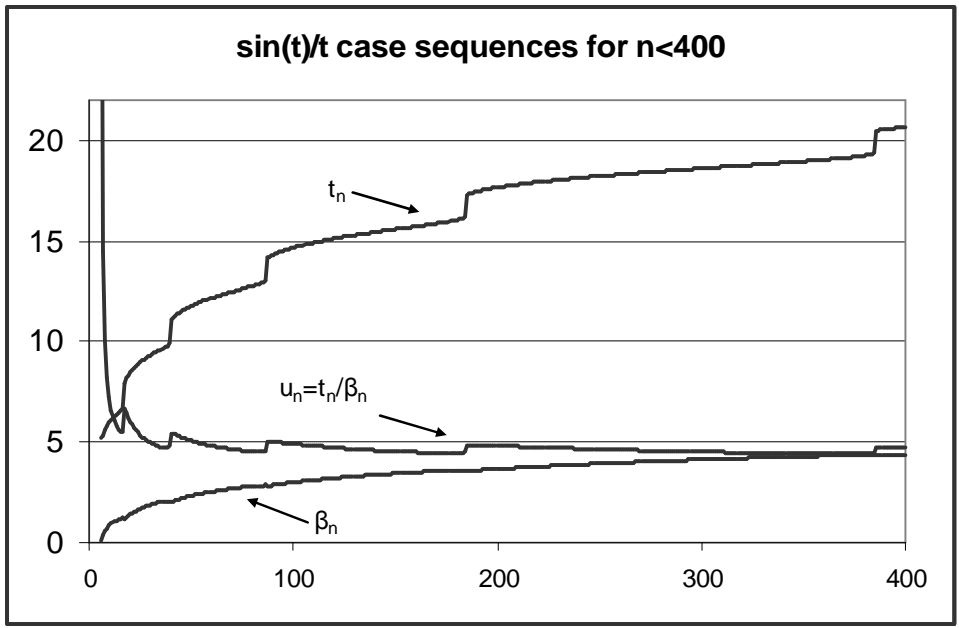

Figure 5.1

The minimal beta-sequence is seen to be rather smooth and if we try to fit the above graph already in the small sub-range up to n = 100, then we find the following fit to the onset of the real root regime

$$b_n \approx 1.02(\log(n+1))^{0.99} - 1.60 \qquad (5.2)$$

which is a first suggestion that the beta sequence might actually have a simple one times log(n) behaviour. The t-graph can of course not be fitted as well by such a simple formula because of the jumps, but overall it looks as if $t_n$ might also be proportional to log(n). Formulated in a different way, the ratio $u_n = t_n / b_n$ seems to approach a non-zero constant, between 4 and 5 (actually it turns out to be ~ 4.4). Numerically this logarithmic behaviour is in fact seen to persist up to very high n, with Fig. 5.2 showing $b_n$ up to n = 2 million, and with Fig. 5.3 showing $t_n$ up to n = 2 million.

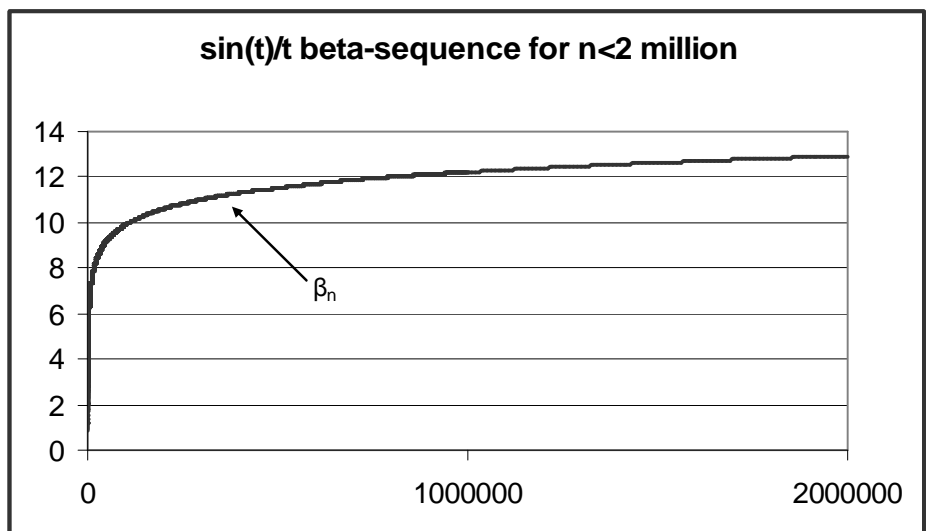

Figure 5.2

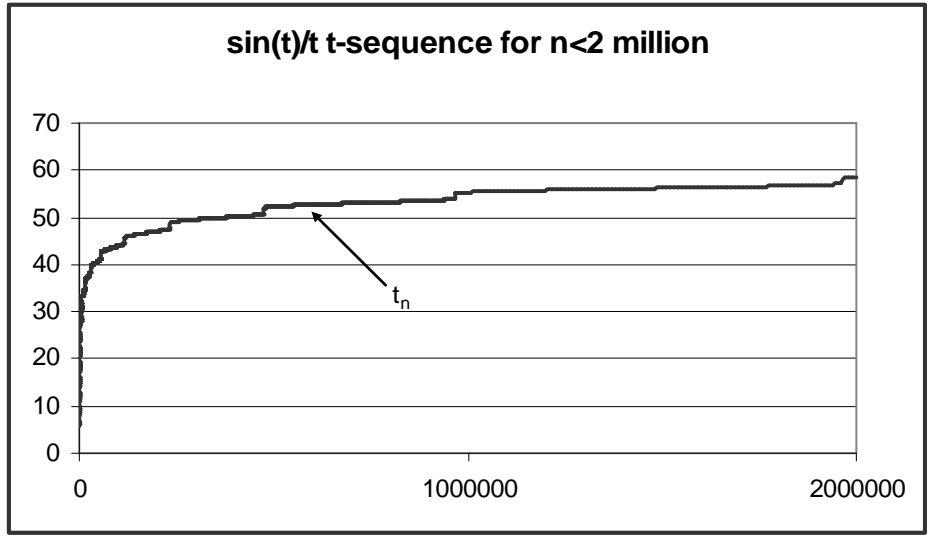

Figure 5.3

The model fit for $b_n$ in the n-range 1000-2 million is

$$b_n \approx 1.0037(\log(n+1))^{0.9996} - 1.61 \tag{5.3}$$

which is very close indeed to the expression above for smaller n.

It is very instructive to compare certain indicators of the general beta and t-increment equations, as well as of the asymptotic form of these equations. In this respect, we may start by noting that the numerator and the denominator terms of the right hand side of the beta-increment equation (in the $Q_{n+1}$ term) have opposite signs, except perhaps in the vicinity of the zeros (or numerically speaking, close to zeros) of these terms, so that the beta-increment is mostly positive, and therefore the $b_n$ sequence is overall increasing. The $t_n$ sequence is strictly increasing, which is also due to the fact that the numerator and the denominator terms of the right hand side of the t-increment equations have opposite signs. In both cases, the key signature for a special sequence behaviour is due to the Pochhammer polynomial $P_{n+1}^{+}(t_n/b_n)$ in the denominator term $\Xi_{n+1}(t_n, b_n)$.

For example, let us examine the graph of $\Xi_4(t, \beta_4)$ where $b_4 \sim 0.074$ in the neighbourhood of the double root $t_4 \sim 3.46$ (see Fig. 5.4). For $b = b_4$, then by construction we have a real double root, but if we increase $b$ slightly to $b = 0.08$, then the $\Xi_4(t, b)$ graph is lowered and the real double root splits into two distinct real roots. On the other hand, if we decrease $b$ slightly to $b = 0.07$, then the $\Xi_4(t, b)$ graph is raised and the real double root splits into a pair of complex roots.

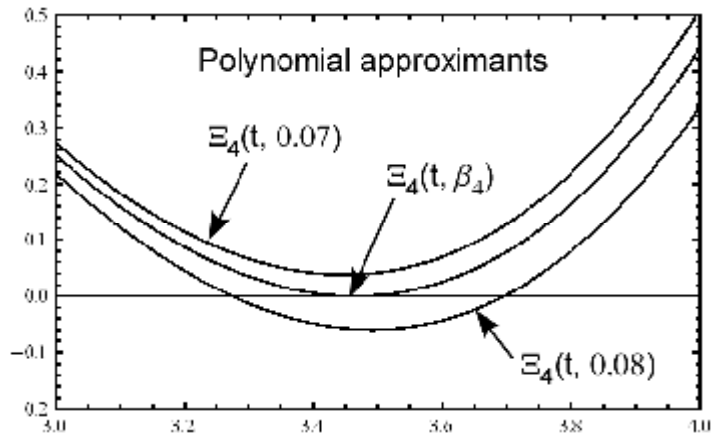

Figure 5.4

As one increments n from 4 to 5, keeping $b_4$ fixed, then one notes that $\Xi_5(t, b_4)$ moves further up above the t-axis (see Fig. 5.5) because the incremental term $P_5^+(t/b_4)$ is positive in this neighbourhood. The iterative mechanism of the beta and t-increment equations is simply to first increase beta so as to push the $\Xi$-graph closer to the t-axis, i.e. making $\Xi(t)$ smaller at $t_4$, and secondly to increase t so as to make $\Xi'(t)$ smaller. This process is iterated until the double root equations are fulfilled at the desired precision. The double root here corresponds to a local minimum of the polynomial.

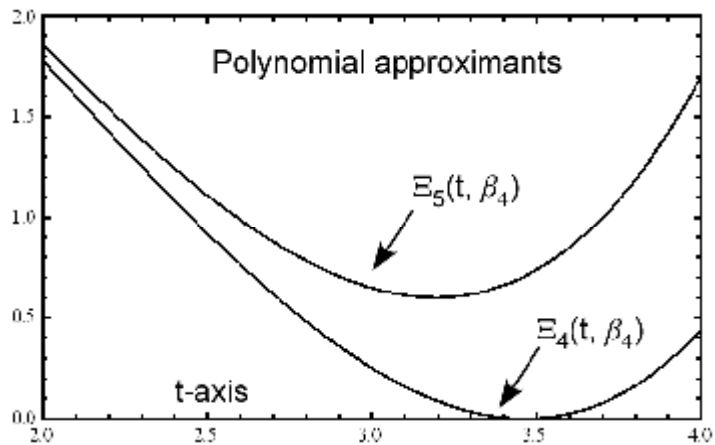

Figure 5.5

The graph of $\Xi_5(t, b_5)$ is shown in Fig. 5.6. However, when we now increment n from 5 to 6, then one observes that the graph of $\Xi_6(t, b_5)$ moves slightly below the axis. This happen because the new incremental term $P_6^+(t/b_5)$ is negative in the double root neighbourhood, in contrast to being positive for lower n. The change of sign of the Pochhammer polynomial is therefore a signature of something special happening. In fact, the polynomial $\Xi_6(t, b_5)$ has a local maximum around t = 5 and when beta is changed, actually lowered somewhat from $b_5 \sim 0.268$ to $b_6 \sim 0.126$, then the $\Xi_6(t, b_5)$ graph is lowered towards the t-axis and the double root jumps from $t_5 \sim 3.47$ to $t_6 \sim 5.17$, which now corresponds to a local maximum of the polynomial $\Xi_6(t, b_6)$.

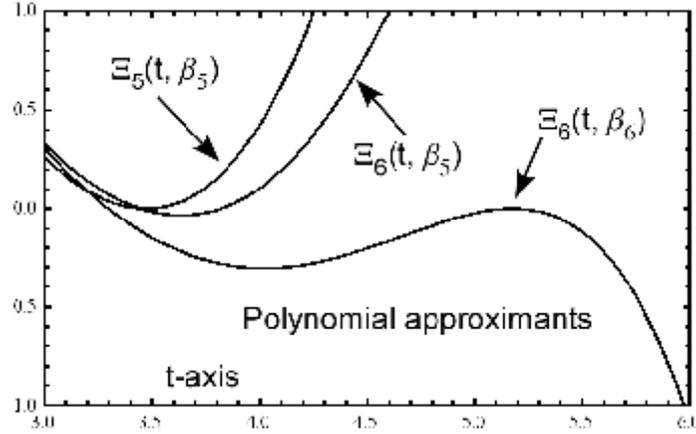

Figure 5.6

This process repeats itself ad infinitum. One notes that $b_n$ and $t_n$ continue to increase smoothly with n until the incremental Pochhammer polynomial term comes close to zero and changes sign. This signals the imminence of a t-jump, which asymptotically is approximately 1.12, and the change-over from a polynomial local maximum to a local minimum, or vice versa.

With the knowledge of the precise values of the minimal beta-sequence and its corresponding t-sequence, for example as found in the range of n up to 400, it is instructive to investigate graphically what happens to the different terms of the right hand side of the asymptotic form of the beta-increment equation, which may be rewritten in the explicit form

$$\Delta b_n = -B_{n+1}(b_n) \frac{P^+_{n+1}(t_n / b_n)}{b_n \frac{\partial}{\partial b} \Xi_{n+1}(t_n, b_n)} \quad (5.4)$$

where

$$B_n(b) = \frac{4}{n+1}(1-e^{-b})^{n+1} \quad (5.5)$$

The overall $B_{n+1}(b_n)$ factor is positive and becomes smaller for higher n. The Pochhammer polynomial factor $P^+_{n+1}(t_n / b_n)$ should asymptotically have an oscillatory behaviour like $\cos(u_n \log(n) + j\,(u_n))$. From the Fig.5.7, with the data points shown as small icons, this is actually seen to be the case. The same oscillatory behaviour is observed for the denominator term which is proportional to a first beta-derivative of the polynomial $\Xi_{n+1}(t_n, b_n)$. As already anticipated, we see that numerator Pochhammer term and denominator beta-derivative term have an opposite oscillatory behaviour in accordance with the beta-sequence being mostly increasing. These two oscillatory terms are quite well synchronized, being close to zero at the same t values, namely the t-jump values.

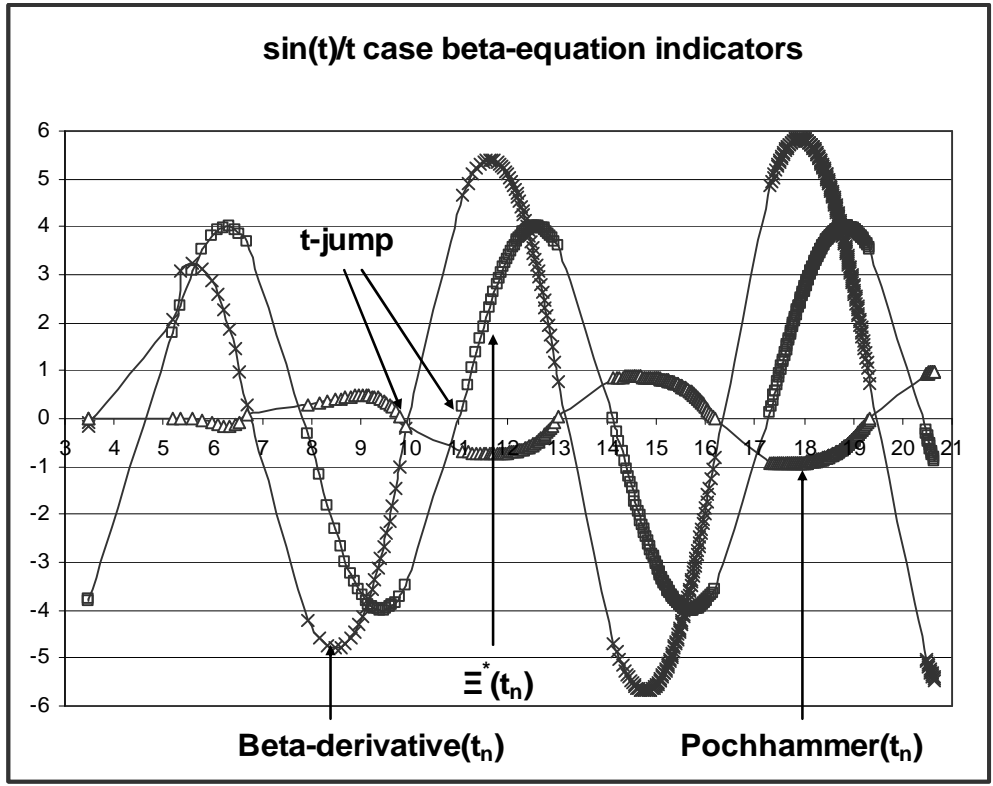

Figure 5.7

Moreover it was noted earlier that the beta-derivative term is expected to be related asymptotically to the function $\Xi^*(t) = \Xi(t) + t\Xi'(t)$ which, with $\Xi(t) = 4\sin(t)/t$, evaluates to exactly $4\cos(t)$. This expectation is actually confirmed by the graph of $\Xi^*(t_n)$ in Fig. 5.7. More precisely we have

$$b_n \tfrac{\partial}{\partial b} \Xi_{n+1}(t_n, b_n) \approx 4\cos(t_n + t_n) \tag{5.6}$$

where $t_n$ is the t-jump which is seen to approach a constant 1.12 for large n.

Similarly we may investigate the behaviour of the corresponding indicator terms in the t-increment equation, which is a little more complex. For short, let us just note that the key numerator term is the same oscillatory Pochhammer polynomial factor, as in the beta-increment equation, and that the denominator term is a second derivative $\Xi''_{n+1}(t_n, b_n)$. The sign of the latter is an indicator of whether we are at local maximum or minimum of the polynomial, and this was precisely seen above to be synchronized with the change of sign of the Pochhammer polynomial. In Fig. 5.8 we show the graph of $P^+_{n+1}(t_n/b_n)$ and that of $\Xi''(t_n)$, which as expected show an opposite oscillatory behaviour, and with a synchronization of the zeros modulo a small t-jump shift.

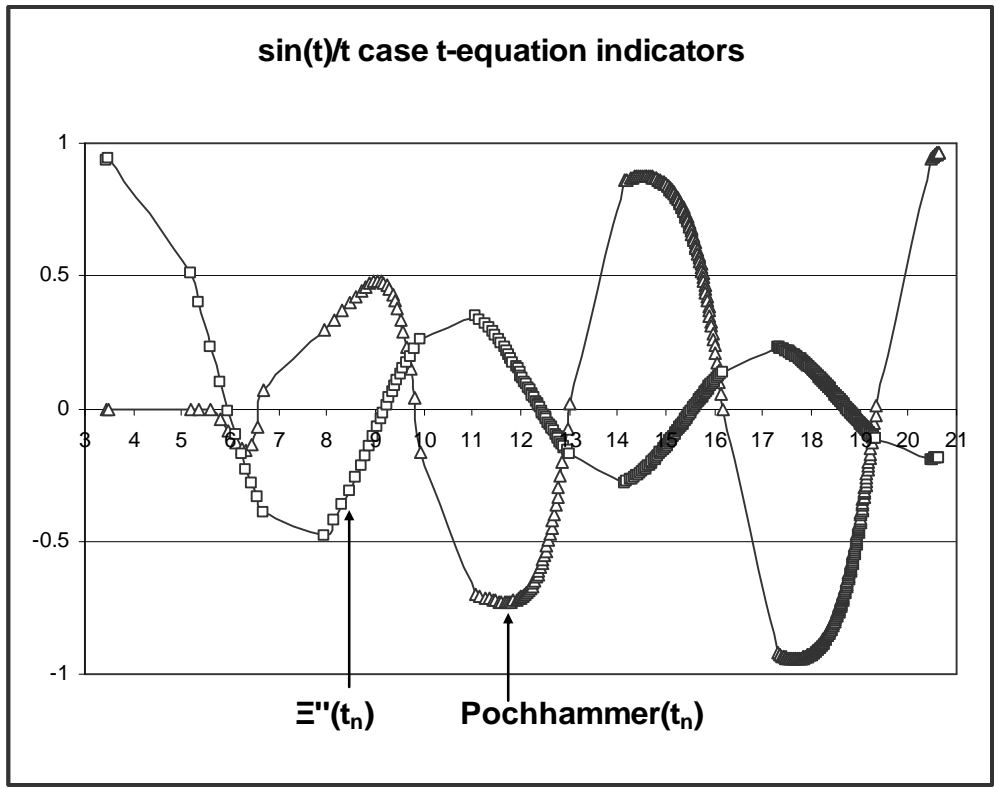

Figure 5.8

The above analysis of the key indicators of the asymptotic beta and t-increment equations was conducted in the relatively low n range up to 400. As n increases we expect these first order equations to become better and better, so the observed graphical features of the numerator and denominator synchronization are also expected to become better and smoother. We will not here discuss the effect of the overall factor $B_{n+1}(b_n)$, as well as the residual constants of the numerator and denominator cancellation mechanism, since this can only be done properly within a purely formal analysis [see Ref. 5].

It is not very different to carry through the numerical analysis of other cases of $A_I(x)$ having compact support, than the above case of $A_I(x) = 1$ for x in the interval [1,exp(2)] corresponding to $\Xi(t)$ = 4sin(t)/t. For example, if one chooses $A_I(x) = \cos(\log(x)p/4)$ on the same interval, then one obtains

$$\Xi(t) \;=\; \frac{2p\cos(t)}{(p/2)^2 - t^2} \tag{5.7}$$

which decreases quadratically with t and has an infinite number of simple zeros, as in the sin(t)/t case. Fig. 5.9 shows graphically the beta and t-sequences in this case for n up to 300. The parameter values are of course somewhat different compared to the $A_I(x)=1$ case, but overall the sequence behaviour is quite similar, e.g. concerning t-jumps.

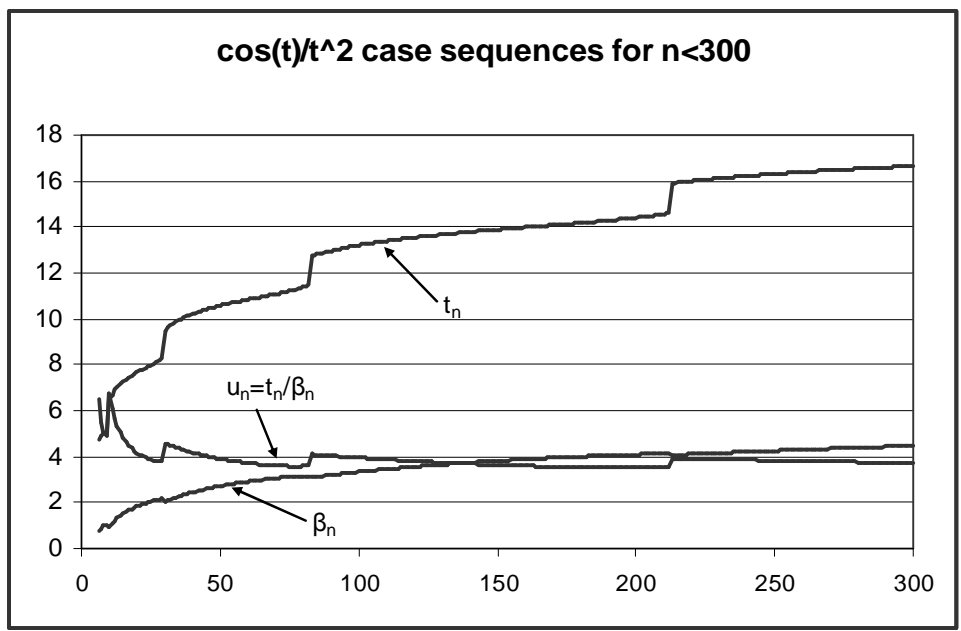

Figure 5.9

## 6. Reference cases with exponential decrease

We have observed a quite rich structure of the beta and t-sequences for some representative cases when $A_I(x)$ has compact support, and one may naturally ask the question whether the same will be the case when $A_I(x)$ has a genuine exponential decrease. Superficially this appear not to be the case when examining some reference cases which are not too dissimilar to the Riemann $A_I(x)$, for example the case $A_I(x) = \exp(-a(x+1/x))$ with a=1 which leads to $\Xi(t)$ being a Bessel function $2K_{it/2}(2)$. In fact, the Riemann $A_I(x)$ is of exponential order 1 and type $\pi$, while the Bessel $A_I(x)$ here is of order 1 and type 1, so some similarity of beta and t-sequences might reasonably be expected.

For this Bessel case, we may start by evaluating $b_{10} \sim 2.12$ and $t_{10} \sim 7.07$ and, proceeding iteratively as before, we find the following graph in Fig. 6.1 for the $b_n$, $t_n$ and $u_n = t_n / b_n$ sequences for n<200:

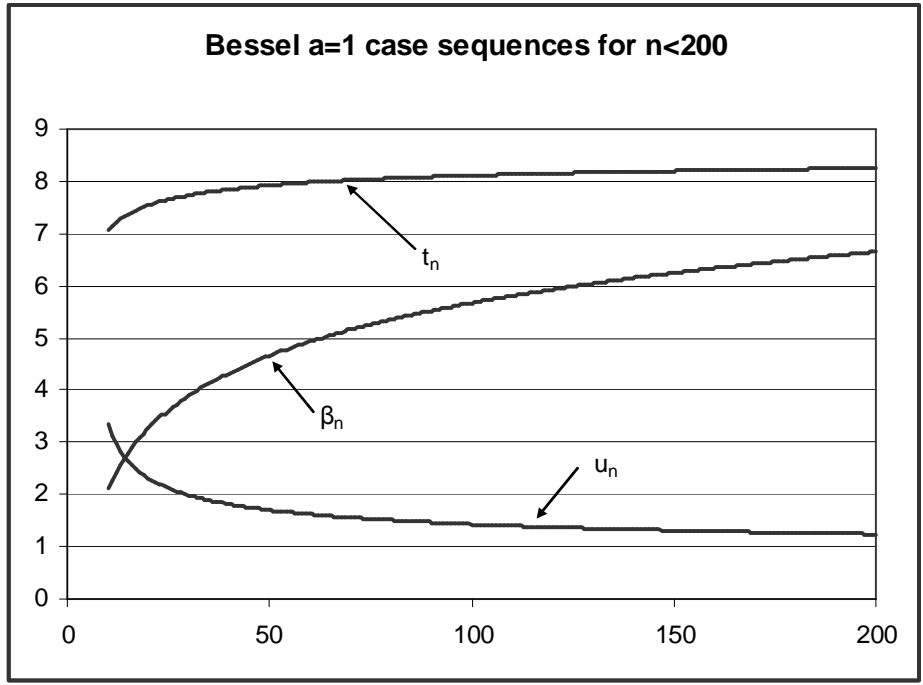

Figure 6.1

We note that $b_n$ increases steadily, while $t_n$ increases much more slowly, and that the $b_n$ ratio decreases, but so far without any clear indication about whether $u_n$ might tend to a non-zero constant, like in the sin(t)/t case, or perhaps towards zero. The interesting feature of the $b_n$ sequence is that it may be well fitted in the studied n range<200 by a logarithmic power function of the form

$$b_n \approx 3.95(\log(n+1))^{0.63} - 4.74 \tag{6.1}$$

It thus appears that there is no t-jump structure, like there was in the sin(t)/t case, and that the $b_n$ growth rate is sub-logarithmic, in contrast to the sin(t)/t case where it appeared to be logarithmic. Of course, this observation can only be indicative and alternative fits are possible in this relatively low n range, but the significant feature is less the particular form of the functional fit than its apparent sub-logarithmic growth rate. Let us also note that the $t_n$ sequence seems to increase much slower than the $b_n$ sequence with a close to logarithmic relationship, and that moreover in this n-range all $t_n$ are confined well below the lowest zero of the Bessel $\Xi(t)$ located at $t \sim 8.81$. This is significant, when recalling from the sin(t)/t case, that t-jumps appeared to happen in between the $\Xi(t)$ zeros.

It turns out that pushing the computer calculations to very large n of the order of a million still does not reveal any t-jump structure. The $b_n$ sequence continues to increase smoothly as shown in the graph of Fig. 6.2, as does the $t_n$ sequence as shown in the graph of Fig. 6.3.

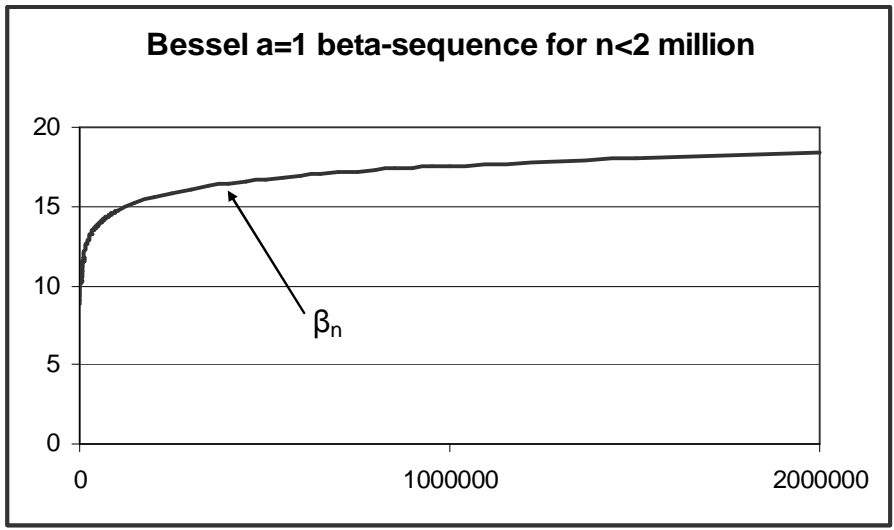

Figure 6.2

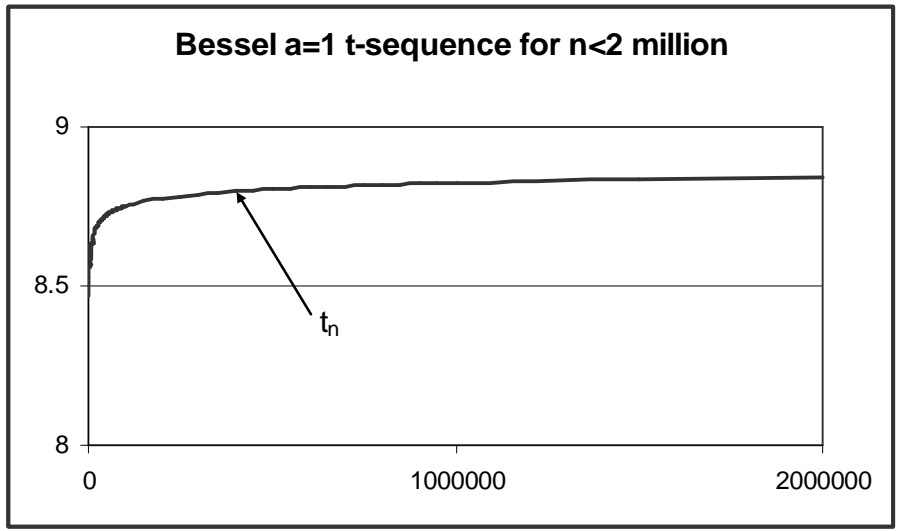

Figure 6.3

The best logarithmic power fit in the n range 1000-2 million for the beta-sequence is

$$b_n \approx 1.72(\log(n+1))^{0.908} - 1.14 \qquad (6.2)$$

which is quite different from the parametrization above for low n. But significantly, the overall growth rate of the $b_n$ sequence is still sub-logarithmic. There are many possibilities of parametrizing a sub-logarithmic behaviour of an increasing sequence $l_n$, but it turns out that it is valuable for the formal analysis of the minimal beta sequences [see Ref. 4] to examine an implicit $l_n$ definition of the general form

$$\log(n) = l_n \, g(l_n) \qquad (6.3)$$

where g is some slowly increasing function. For example, if $g(l) = \log(l)$ and $\log(\log(l))$, then we get the useful sub-logarithmic functions [see Ref. 5] $l_n$ = sublog(n) and sublogxl(n), respectively. These functions fulfill the growth rate relation sublog(n)<<sublogxl(n)<<log(n). It is not difficult to see that by enlarging this definition of the function g with additional constants and sub-dominant terms, then it is possible to propose several alternative parametrizations of the sequence data presented in the graphs above, other than the logarithmic power ones considered earlier.

One might speculate that the absence of any particular structure of the above beta and t-sequences might be due the special exponential form of $\Xi(t)$ and that consequently some sequence features might possibly somehow only emerge very slowly on a logarithmic scale of the relevant parameters. We will therefore use the parameter freedom inherent in the Bessel reference case to choose a different parameter in $\Xi(t) = \exp(-a(x+1/x))$ which might show a more interesting sequence structure for low n. More generally we have $\Xi(t) = 2K_{it/2}(2a)$, and it turns out that when choosing a to be much smaller than 1, say a=0.005, then something interesting indeed happens (see Fig. 6.4). Starting at n=4 we find $b_4 \sim 0.044$ and $t_4 \sim 1.34$ which is already located between the first and second zero of this $\Xi(t)$, at t ~ 1.29 and 2.47, respectively. Then following some small increases of β and t, at n=7 we see t jump from 1.43 to1.92 at n=8. Subsequently there are t-jumps at n = 32 and 161, as well as further on.

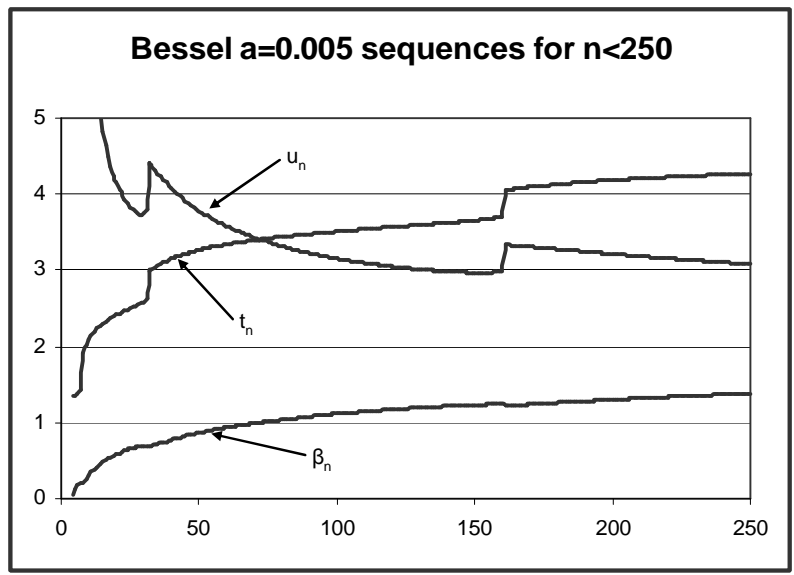

Figure 6.4

In this Bessel $\Xi(t)$ case of a=0.005, we recover all the peculiar sequence features of the sin(t)/t case. This $\Xi(t)$ also has infinitely many zeros and the t-jumps of the polynomial approximants happen in between these zeros, so there is good reason to anticipate that the t-jump process will actually continue ad infinitum. But we also note that the values of n at which there are t-jumps increase in a multiplicative manner. In the sin(t)/t case there was an approximate factor of 2 between successive jumps, but in the Bessel case, this factor is bigger, and in fact it seems to become bigger and bigger. This circumstance may of course eventually pose a practical problem for computer calculations since the identification of each successive jump feature will require steeply increasing computing resources.

The above numerical observations lead us to assert that, in spite of the fact that Bessel a=1 case did not reveal any t-jump structure as far as n=2'000'000, this $\Xi(t)$ case should eventually also show similar features as the other reference cases. However, it requires more analysis to predict at exactly which higher n one will start seeing t-jump features happening.

Even for low n, the Bessel a=0.005 case does reveal some crucial information about the asymptotic indicators of the fundamental mechanism at work in the beta and t-increment equations. This is clear when plotting in Fig. 6.5 the same indicators, i.e. beta-derivative, Pochhammer and $\Xi^*(t)$ terms as done in the sin(t)/t case:

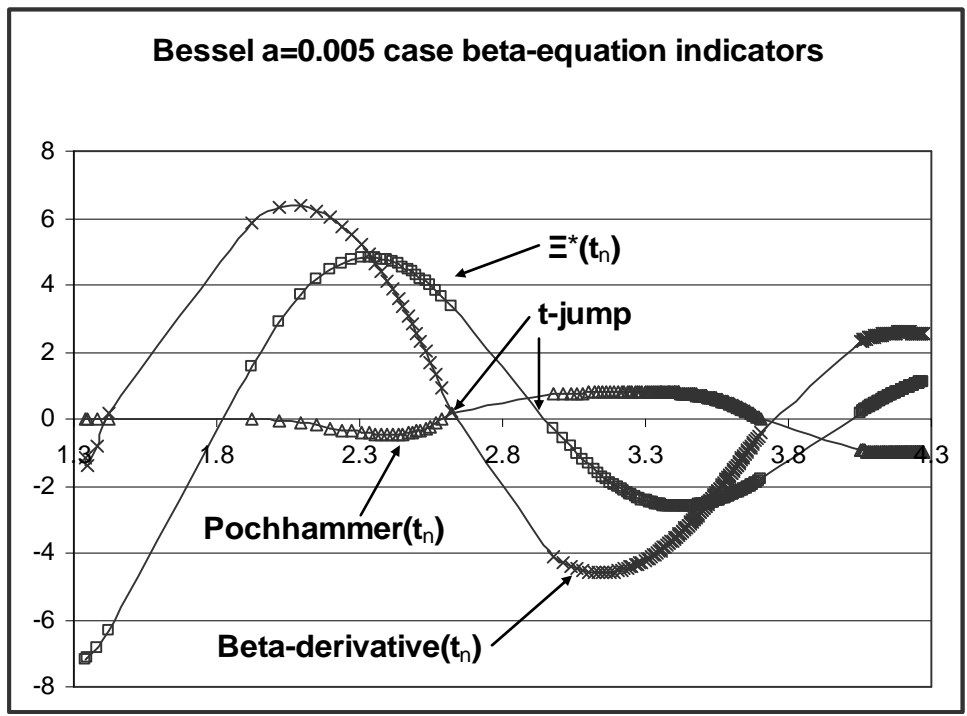

Figure 6.5

We again note that the numerator term on the right hand side of the beta-increment equation is a Pochhammer polynomial which is oscillatory and opposite to the beta-derivative term of the denominator. The zeros of these terms are synchronized, and moreover the beta-derivative term is related to the $\Xi^*(t)$ quantity modulo a phase shift $\tau_n$ which is precisely the t-jump. In contrast to the sin(t)/t case where $t_n$ was seen to approach a constant, the t-jumps here become closer and closer to zero for large n. This can be understood in terms of the asymptotic behaviour of $\Xi(t)$, which is of an exponentially damped oscillatory form, and the fact that the number of zeros of $\Xi(t)$ grows faster than linearly. It is thus pretty safe to infer that this kind of asymptotic mechanism will also be operational for the Bessel case a=1, even if it has not yet been observed computationally.

The other exponential reference cases concern the simpler choice of $A_I(x) = \exp(-ax)$, say with a=1 leading to $\Xi(t)$ being a symmetrized incomplete gamma function $\Xi(t) = \Gamma(it/2,1) + \Gamma(-it/2,1)$. Let us recall that this function is positive and only has a real zero at infinity. This of course does not prevent the polynomial approximants from having roots and the standard procedure from being applicable to finding the beta and t-sequences as above. Starting with $b_{10} \sim 2.76$ and $t_{10} \sim 8.32$, one may proceed further to note that the $b_n$ sequence has a smooth sub-logarithmic behaviour, while the $t_n$ sequence appears to increase much more slowly as shown in the Fig. 6.6

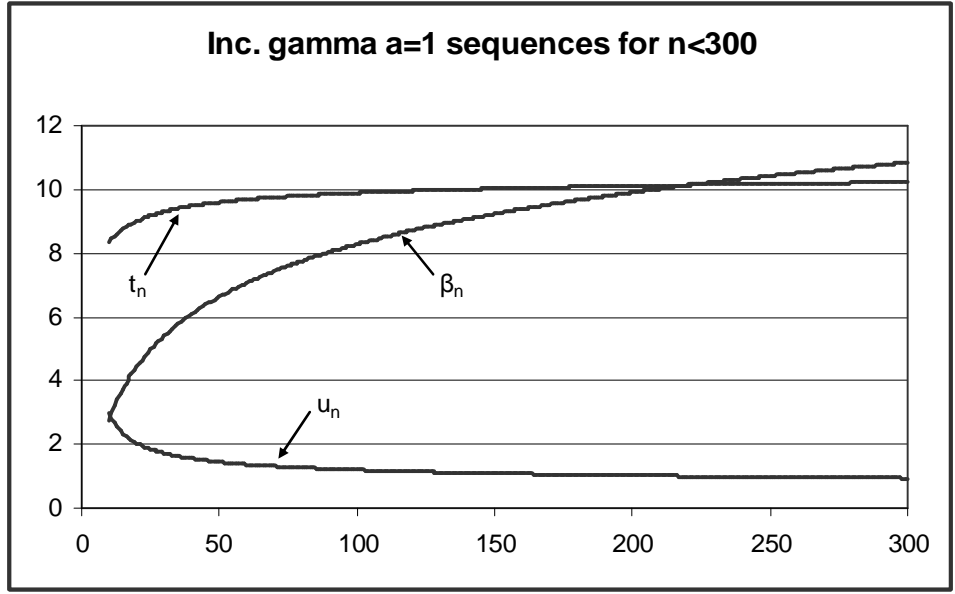

Figure 6.6

The best fit to the beta-sequence in this n-range up to 300 is given by

$$b_n \approx 3.89(\log(n+1))^{0.865} - 4.45 \qquad (6.4)$$

It appears to have a higher growth rate than that of the Bessel a=1 case, but it still has a distinct sub-logarithmic behaviour. The $t_n$ sequence is growing very slowly, having no t-jumps, and seems to be confined below a certain value like in the Bessel a=1 case, even if in this case, $\Xi(t)$ does not have any zero. The $u_n$ ratio also seems to approach zero.

As in the Bessel case, the smooth beta and t-sequence situation changes radically if instead of a=1 we choose a much lower value, say a=0.01. In this case, the corresponding $\Xi(t)$ has exactly 8 zeros located at 1.49, 2.81, 4.04,…Starting at n=4 we find $b_4 \sim 0.054$ and $t_4 \sim 1.53$ which is already located between the first and second zero of $\Xi(t)$. Then following some small increases of $b$ and t, at n=7 we see t jump from 1.62 to 2.12 at n=8 (see Fig. 6.7). Subsequently there are t-jumps at n = 36 and 248, as well as further on.

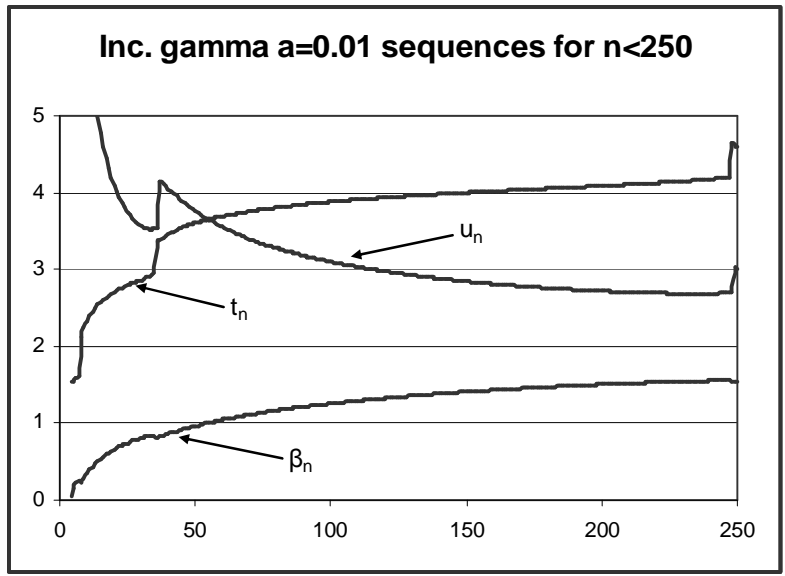

Figure 6.7

Thus we recover the same t-jump features as seen in the sin(t)/t and Bessel a=0,005 cases, indicating that the same fundamental mechanism of the beta and t-increment equations is operational. However, here this situation cannot persist in exactly the same way since $\Xi(t)$ only has a finite number of zeros. After having observed as many t-jumps as there are zeros of $\Xi(t)$, then the underlying mechanism must change character [see Ref. 4] and will resemble more the situation of the incomplete gamma function case with a=1 where $\Xi(t)$ had no zeros. But until it does so, then we find exactly the same behaviour of the beta-sequence indicators (see Fig.6.8) as previously found:

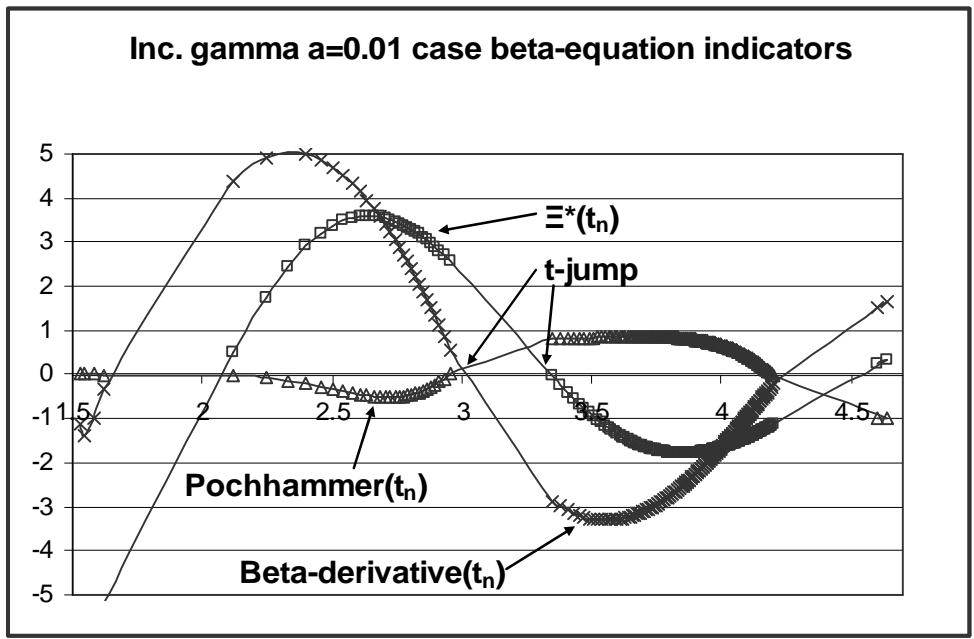

Figure 6.8

## 7. The Riemann case

The exponential reference cases are useful for acquiring a preliminary understanding of what may be the nature of the beta and t-sequences for the Riemann case where $A_I(x)$ is expressed in terms of the elliptic theta function as an infinite sum of exponential terms with power factors. Retaining only the leading exponential behaviour of the Riemann $A_I(x)$, and substituting x→x+1/x to assure that $\Xi(t)$ has an infinite number of zeros, then one obtains the Bessel function type reference case. On the other hand, retaining only a finite number of terms in the Riemann $A_I(x)$ series, then one recovers the reference case of a sum of incomplete gamma functions having only a finite number of zeros.

The numerical analysis of the beta and t-sequences of the Riemann case can of course in practice, i.e. for finite n, be conducted by keeping only a finite number of terms of the Riemann $A_I(x)$ series. So as usual we may start by evaluating $b_{10}$ ~ 3.23 and $t_{10}$ ~ 10.63 and proceed iteratively as done in the other reference cases to find the low n (up to 200) behaviour of the sequence as shown in Fig.7.1.

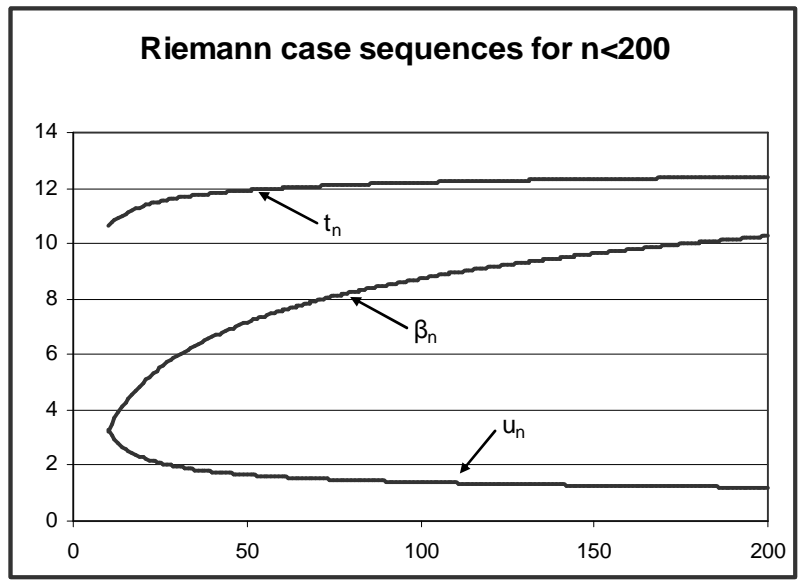

Figure 7.1

Like in the Bessel a=1 case, we note a smooth increasing behaviour of the $b_n$ sequence, a much slower growth of the $t_n$ sequence, and the $u_n$ sequence seemingly approaching zero. There are no t-jumps and in the studied n-range, the $t_n$ appear to be confined below the lowest zero of $\Xi(t)$ at t~14.1. The best fit using a logarithmic power parametrization is

$$b_n \approx 5.58(\log(n+1))^{0.66} - 6.76 \qquad (7.1)$$

which is seen to be distinctly sub-logarithmic much like in the Bessel a=1 case.

It turns out that pushing the computer calculations to very large n of the order of a million still does not reveal any t-jump structure. The $b_n$ sequence continues to increase smoothly as shown in the graph of Fig. 7.2, and there are no t-jumps, as shown in Fig. 7.3.

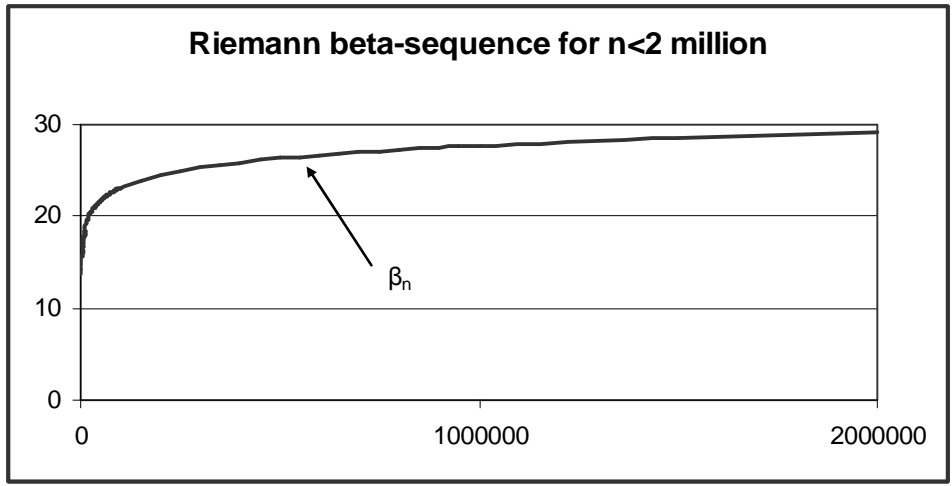

Figure 7.2

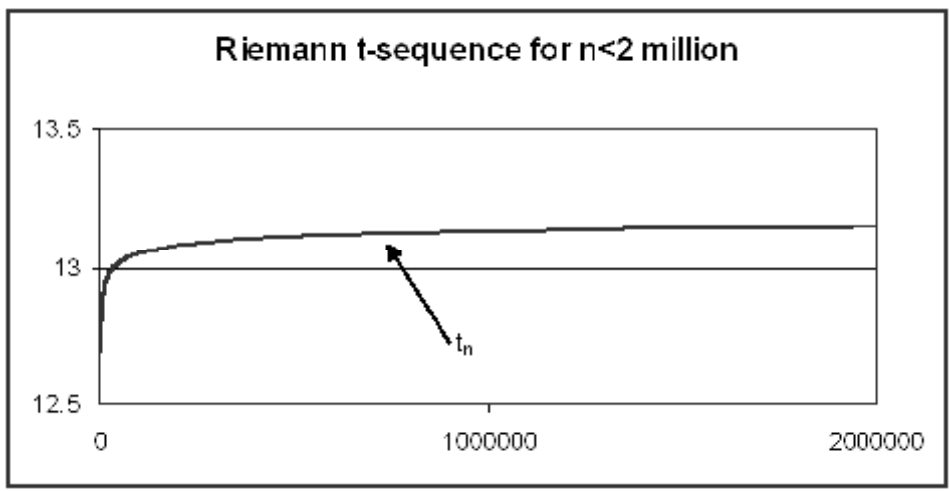

Figure 7.3

The best logarithmic power fit in the n range 1000-2 million for the beta-sequence is

$$b_n \approx 2.57(\log(n+1))^{0.928} - 1.75 \qquad (7.2)$$

which is quite different from the parametrization above. But significantly the overall growth rate of the $b_n$ sequence is still sub-logarithmic, like in the Bessel a=1 case. As note earlier, it is not difficult to envisage different alternative sub-logarithmic fits to the beta-sequence data, for example of the sublog(n) and sublogxl(n) type discussed. However, to do so it would be useful to justify one form or the other by arguments from a formal analysis [see Ref.4], but this goes beyond the scope of the present paper.

In contrast to the exponential reference cases, there is no free parameter to change in the Riemann $A_I(x)$ case because its exponential type is fixed to a = $p$. Nevertheless, from the lessons of the Bessel and incomplete gamma function cases, it is pretty safe to assume that a similar t-jump structure will eventually be seen at some point above n = 2 million, quite possibly very much above this. Moreover, the same fundamental mechanism involving the terms of the beta and t-increment equations can also be assumed to be operational for the Riemann case.

The open question is of course whether all of this may lead to the conclusion that the minimal beta-sequence of the Riemann case asymptotically has a truly sub-logarithmic growth [see Ref. 4 for a discussion of this point]. The main difficulty in trying to come to such a conclusion is clearly that our present knowledge about the Riemann $\Xi(t)$ function is incomplete in what concerns its asymptotic behaviour for large t, as well as its distribution of real zeros. These two features were seen to be quite important for the understanding of the asymptotic interplay of the indicator terms in the beta-increment equation.

## 8. Exceptional cases

The family of entire functions characterized by a certain A(x) is quite broad and some cleverly engineered choices of A(x) can lead to exceptional properties of the corresponding $\Xi(t)$ and beta-sequences which are somewhat surprising. In particular, this is so for the $A^{(1)}(x,k)$ and $A^{(2)}(x,k)$ cases, mentioned earlier, which are interesting because their corresponding $\Xi(t)$ have complex zeros off the critical line. However, our current understanding of these $\Xi$-functions is very far from being complete, specifically concerning the distribution of their real zeros and the asymptotic behaviour for large real t. This lack of understanding also implies that it is not so easy to fully anticipate the behaviour of the corresponding beta and t-sequences.

The $A^{(1)}(x,k)$ case for k=5, for example, related to the Ramanujan tau function has some points in common with the sin(t)/t case. In the range of n<200, we may start by calculating $b_{10} \sim 1.006$ and $t_{10} \sim 3.47$ and proceed to find that the $b_n$ sequence is increasing smoothly, except for n = 66 where it decreases somewhat, as seen in Fig. 8.1. At this point, the $t_n$ root is seen to jump upwards within the interval of the lowest $\Xi(t)$ zeros located at 3.83 and 6.41

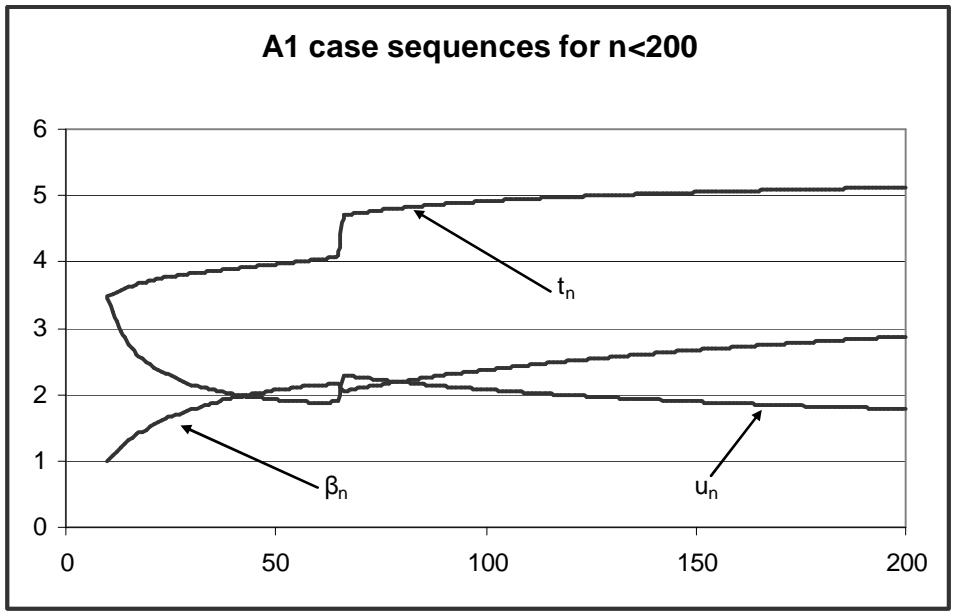

Figure 8.1

The best fit to the growth rate of the $b_n$ sequence shows the following supra-logarithmic behaviour

$$b_n \approx 0.302(\log(n+1))^{1.295} + 0.213 \qquad (8.1)$$

The growth rate of the $t_n$ sequence appears to be not too different from that of the beta-sequence.

For the $A^{(2)}(x,k)$ case for k=5, for example, with no real zeros as mentioned earlier, we may start by calculating $b_{10} \sim 2.11$ and $t_{10} \sim 7.31$ and proceed to find and a rather smooth monotonous behavior of the $b_n$ sequence and a $t_n$ sequence with no t-jumps as seen in Fig. 8.2

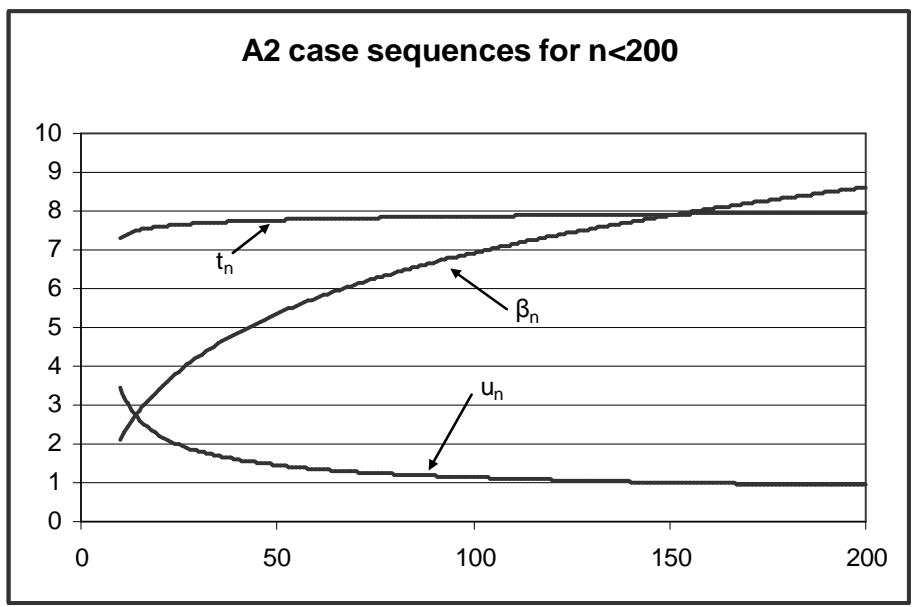

Figure 8.2

The best fit to the growth rate of the $b_n$ sequence has the following supra-logarithmic behaviour

$$b_n \approx 1.099(\log(n+1))^{1.321} - 1.364 \qquad (8.2)$$

The growth rate of the t-sequence appears to be much slower than that of the beta-sequence.

## 9. Conclusions

The numerical analysis of different cases of admissible A(x) functions defining the broad family of entire functions under study has revealed a quite interesting structure of the associated minimal beta and t-sequences. The real double root condition defining these sequences gives rise to two equations for the beta and t-increments which are the basis for both a computerized numerical analysis and a more formal asymptotic analysis.

The numerical analysis confirms the generic feature of $b_n$ sequences which are globally increasing, apart from small patches of slight decrease, as well as globally increasing $t_n$ sequences characterized by more or less big t-jumps. In some cases, these particular features are seen already for very small n, as well as for higher n, and in some cases, including the Riemann A(x) case, they are not seen at all in an n-range up to 2 million. However, there is quite good reason to believe that similar sequence features are present for all A(x) cases, even if the features may sometimes be manifest only for very large n, and that the features are due the same underlying mechanism of the asymptotic increment equations.

For the purpose of the present paper, the most interesting feature of the sequences is the growth rate of the minimal $b_n$ sequence. The application of the Hurwitz theorem of complex analysis is possible if the polynomial approximants $\Xi_n(t, b_n)$ actually converge to $\Xi(t)$, as the fixed $b$ approximants $\Xi_n(t, b)$ do. It is outside the scope of the present paper to discuss the details of this convergence issue [see Ref. 4], but let us just state that there are convincing indications that the $b_n$ growth rate of log(n) represents a dividing line under which there is convergence, and above which there is non-convergence of the $\Xi_n(t, b_n)$ approximants.

For the case of A(x) with compact support, like the $\Xi(t) = 4\sin(t)/t$ case, the situation is very special and the analysis is rather clear. The numerical analysis has provided strong support for a minimal $\beta_n$ growth rate of exactly log(n). A formal asymptotic analysis [see Ref. 4] substantiates this observation and makes it even more precise, so as to state that $b_{\min,n} = \log(n/n_0)$ where $n_0 \sim 5.1$. It can also be shown that there is convergence of the $\Xi_n(t, b_n)$ approximants if $b_n$ grows at most like log(n) + const. Therefore the Hurwitz theorem can be applied to this case albeit, as already noted, the subsequent conclusion about the real zeros of $\Xi(t)$ does not provide any new information.

The cases of A(x) with genuine exponential decrease are more complex. For the exponential reference cases, we have seen numerical indications for both sub-logarithmic $b_n$ growth rate, i.e. for the Bessel, incomplete gamma, and Riemann cases, as well as for supra-logarithmic $b_n$ growth rate, i.e. for the exceptional cases of Dirichlet type. But numerical indications for the large n of the $b_n$ sequences do of course not amount to a statement about the true asymptotic behaviour. Nevertheless, further numerical analysis involving very extensive computer resources may be useful for adding more substance to the indications discussed in the present paper.

A very significant conclusion of the present numerical analysis concerns the information provided about the indicator terms in the asymptotic beta and t-increment equations. Already for low n, this information allows to identify a fundamental mechanism of opposite oscillatory behaviours and synchronized zeros, which can be readily incorporated into a formal asymptotic analysis [see Ref. 4]. The manifestation of this mechanism is most transparent for the cases of $\Xi(t)$ with an infinite number of zeros, like the sin(t)/t, Bessel and Riemann case. For the cases of $\Xi(t)$ with a finite number of zeros, or no zeros, like the incomplete gamma and the truncated Riemann cases, the mechanism is a little different [see Ref. 4], which may also be seen to explain why the numerical results for the $b_n$ growth rate for these cases, while still appearing to be sub-logarithmic, are distinguishable from the other exponential cases.

In summary, the numerical analysis of minimal beta-sequences has been seen to reveal a very rich structure of features with potentially important applications, both for computer analysis and for formal analysis. The impetus for this analysis was due to the existence of a Pochhammer polynomial expansion depending of a continuous parameter $b$, so that one may define $b$-dependent polynomial approximants, the properties of which may be suitably adjusted. Here we studied the possibility of adjusting the polynomial roots to be real, so as to eventually apply the Hurwitz theorem for cases where the approximants can be shown to converge.

It is possible to generalize this approach in different ways. Fundamentally the Pochhammer polynomial expansion originated in a polynomial Euler formula [see Ref. 4] and therefore the approximants may be seen to implement a kind of polynomial Fourier analysis. Different types of transforms (Fourier, Laplace and Mellin transforms) may be considered as well as different types of polynomials with appropriate generating function, all dependent on the type of application under study. As shown above, one quite appealing feature of such approaches is that the numerical analysis and the more formal analysis may advantageously be developed in parallel.